\documentclass[12pt]{article}
\usepackage{amsmath,amssymb,amsthm, amscd, pb-diagram}

\begin{document}

\begin{center}
{\bf Toric hyperk\"{a}hler varieties and {\bf Q}-factorial terminalizations}
\end{center}
\vspace{0.4cm}

\begin{center}
{\bf Yoshinori Namikawa}      
\end{center}
\vspace{0.2cm}

{\bf Introduction}. 
\vspace{0.2cm}

A toric hyperk\"{a}hler variety  is defined as a hyperk\"{a}hler reduction of  
a quarternionic space $\mathbf{H}^n$ with the standard hyperk\"{a}hler structure $(g, I, J, K)$, by an action of a compact torus $T^d_{\mathbf R}$ (cf. [BD], [Go], [HS], [Ko 1]).
If we identify $\mathbf{H}^n$ with $\mathbf{C}^{2n}$ by the complex structure $I$, a toric hyperk\"{a}hler variety is regarded as a holomorphic symplectic reduction of $(\mathbf{C}^{2n}, \omega)$ by a Hamiltonian action of an algebraic torus $T^d$. Here $\omega$ is the standard symplectic 2-form on $\mathbf{C}^{2n}$.  Recently it has been studied from various points of view such as symplectic dulality, biratonal geometry and Poisson deformatons ([BLPW], [Nag]). In this article we denote by $Y(A, \alpha)$ a toric hyperk\"{a}hler variety according to [HS].     
Here $A$ is a $d \times n$ integer valued matrix such that $\mathbf{Z}^n \stackrel{A}\to \mathbf{Z}^d$ is a surjection, and $\alpha$ is an element of $\mathrm{Hom}_{alg.gp}(T^d, \mathbf{C}^*)$. Let us introduce an $n \times (n-d)$ integer valued matrix  $B$ 
in such a way that the sequence $$0 \to \mathbf{Z}^{n-d} \stackrel{B}\to \mathbf{Z}^n \stackrel{A}\to \mathbf{Z}^d \to 0 $$ is 
exact. $Y(A, \alpha)$ essentially depends on this $B$ rather than $A$, and certain properties of $B$ often reflect the geometry of $Y(A, \alpha)$. We assume that all row vectors of $B$ are nonzero.
When $\alpha = 0$, $Y(A, \alpha)$ is a conical symplectic 
variety, and $Y(A, \alpha)$ is a crepant, projective partial resolution of $Y(A, 0)$ for any $\alpha$. When one chooses $\alpha$ generic, $Y(A, \alpha)$ 
has only quotient singularities. In particular, if $A$ is unimodular, then $Y(A, \alpha)$ is nonsingular for a generic $\alpha$, and it gives a crepant resolution of $Y(A, 0)$. But, in general, $Y(A, 0)$ does not have any crepant resolution. A substitute for a crepant resolution is a good crepant partial resolution called a {\bf Q}-factorial terminalization. It would be natural to expect that $Y(A, \alpha)$ is a {\bf Q}-factorial 
terminalization of $Y(A, 0)$ if $\alpha$ is chosen generic. But, this is not true because $Y(A, \alpha)$ may possibly have 
singularities of codimension 2 even if $\alpha$ is generic. In this article, we realize $Y(A, 0)$ as another toric hyperk\"{a}hler variety 
$Y(A^{\sharp}, 0)$ so that $Y(A^{\sharp}, \alpha^{\sharp})$ is a {\bf Q}-factorial terminalization 
of $Y(A^{\sharp}, 0)$ for a generic $\alpha^{\sharp}$ (Theorem 11). More precisely, this $A^{\sharp}$ is charactrized by the 
property that $B^{\sharp}$ is the primitivization of $B$. For the notion of a primitivization, see (4.3).    
As an application we give a necessary and sufficient condition for $Y(A, 0)$ to have a crepant projective resolution. 
In fact, $Y(A, 0)$ has a crepant projective resolution if and only if the primitivization $B^{\sharp}$ of $B$ is unimodular (Corollary 13). Moreover, we construct very explicitly the universal Poisson deformation of $Y(A, 0)$ in terms of $A^{\sharp}$ (see (4.1), (4.2) and (4.3)). \vspace{0.2cm}

\S {\bf 1}.
\vspace{0.2cm}

Let $\mathbf{C}^{2n}$ be an affine space with coordinates $z_1$, ..., $z_n$, $w_1$, ..., $w_n$. 
An $n$ dimensional algebraic torus $T^n$ acts on $\mathbf{C}^{2n}$ by 
$$(z_1, ..., z_n, w_1, ..., w_n) \to (t_1z_1, ..., t_nz_n, t_1^{-1}w_1, ..., t_n^{-1}w_n)$$
By an integer valued $d \times n$-matrix $A := (a_{ij})$, we determine a homomorphism of algebraic tori 
$\phi: T^d \to T^n$ by $$ (t_1, ..., t_d) \to (t_1^{a_{11}}\cdot\cdot\cdot t_d^{a_{d1}}, ..., t_1^{a_{1n}}\cdot\cdot\cdot t_d^{a_{dn}}).$$ 
Then $T^d$ acts on $\mathbf{C}^{2n}$ by 
$$(z_1, ..., z_n, w_1, ..., w_n) \to$$ $$(t_1^{a_{11}}\cdot\cdot\cdot t_d^{a_{d1}}z_1, ..., t_1^{a_{1n}}\cdot\cdot\cdot t_d^{a_{dn}}z_n, 
t_1^{-a_{11}}\cdot\cdot\cdot t_d^{-a_{d1}}w_1, ..., t_1^{-a_{1n}}\cdot\cdot\cdot t_d^{-a_{dn}}w_n)$$ 
The homomorphism $\phi$ induces a map of characters: $\phi^*: \mathrm{Hom}_{alg.gp}(T^n, \mathbf{C}^*) \to 
\mathrm{Hom}_{alg.gp}(T^n, \mathbf{C}^*)$. When we identify the character groups respectively with $\mathbf{Z}^n$ and 
$\mathbf{Z}^d$ in a natural way, $\phi^*$ is nothing but the homomorphism $\mathbf{Z}^n \stackrel{A}\to \mathbf{Z}^d$ 
determined by $A$. We assume that the map $A$ is a surjection. Let $B$ be an integer valued $n \times (n-d)$-matrix 
such that the following 
sequence is exact: 
$$0 \to \mathbf{Z}^{n-d} \stackrel{B}\to \mathbf{Z}^n \stackrel{A}\to \mathbf{Z}^d \to 0.$$  
We assume that all row vectors of $B$ are non-zero. 

Define a symplectic 2-form $\omega$ on $\mathbf{C}^{2n}$ by 
$$\omega := \sum_{1 \le i \le n} dw_i \wedge dz_i.$$  Then the $T^d$-action is a Hamiltonian action on $(\mathbf{C}^{2n}, \omega)$. 
Writing $A = (\mathbf{a}_1, ..., \mathbf{a}_n)$ by the column vectors, the moment map $\mu: \mathbf{C}^{2n} \to \mathbf{C}^d$ is 
given by $$(z_1, ..., z_n, w_1, ..., w_n) \to \sum_{1 \le i \le n}\mathbf{a}_iz_iw_i.$$  

{\bf Lemma 1}. {\em $\mu^{-1}(0)$ is of complete intersection. In particular, $\dim \: \mu^{-1}(0) = 2n -d$.} 
\vspace{0.2cm}

{\em Proof}. Put $W := \mathrm{Spec}\mathbf{C}[z_1w_1, ..., z_nw_n]$. Then 
$\mu$ factorizes as $\mathbf{C}^{2n} \stackrel{\iota}\to W \stackrel{\nu}\to \mathbf{C}^d$. Here $\nu$ is determined by the ring homomorphism $$\mathbf{C}[s_1, ..., s_d] \to \mathbf{C}[z_1w_1, ...., z_nw_n], \:\:\ (s_i \to \sum_{j = 1}^n 
z_jw_j a_{ij})$$  and $\iota$ is determined by the inclusion 
$$\mathbf{C}[z_1w_1, ..., z_nw_n] \to \mathbf{C}[z_1, ..., z_n, w_1, ..., w_n].$$
$W$ is an $n$ dimensional affine space, Since $\mathrm{rank}\: A = d$, 
$\nu^{-1}(0)$ is an $n-d$ dimensional linear subspace of $W$. It is easily checked that $\iota$ is a flat map. 
Hence $\mu^{-1}(0) = \iota^{-1}(\nu^{-1}(0))$ has dimension $n + (n-d) = 2n-d$. Since $\mu^{-1}(0)$ is defined by $d$ equations in $\mathbf{C}^{2n}$, we see that $\mu^{-1}(0)$ is of complete intersection. 
$\square$ 
\vspace{0.2cm}

Note that $T^d$ acts on each fiber of $\mu$. 
Put $M := \mathrm{Hom}_{alg.gp}(T^d, \mathbf{C}^*)$. 
For $\alpha \in M$, we define $$X(A, \alpha) := \mathbf{C}^{2n}/\hspace{-0.1cm}/_{\alpha}T^d, \:\:\: Y(A, \alpha) := \mu^{-1}(0)
/\hspace{-0.1cm}/_{\alpha}T^d.$$ 
When $\alpha = 0$, $$ X(A, 0) = \mathrm{Spec}\: \mathbf{C}[z_1, ..., z_n, w_1, ..., w_n]^{T_d}, \:\: Y(A, 0) = \mathrm{Spec}\: \mathbf{C}[\mu^{-1}(0)]^{T_d}.$$ 
Denote by $(\mathbf{C}^{2n})^{\alpha-ss}$ and $(\mu^{-1})^{\alpha-ss}(0)$ respectivelt the $\alpha$-semistable locus for the $T^d$-action on $\mathbf{C}^{2n}$ and $\mu^{-1}(0)$. 
The inclusion maps $(\mathbf{C}^{2n})^{\alpha-ss} \to \mathbf{C}^{2n}$ and $(\mu^{-1})^{\alpha-ss}(0) \to \mu^{-1}(0)$ respectively induce maps $\nu_{X, \alpha}: X(A, \alpha) \to X(A, 0)$ and $\nu_{Y, \alpha}: Y(A, \alpha) \to Y(A, 0)$. They are  
birational projective morphisms. The moment map $\mu: \mathbf{C}^{2n} \to \mathbf{C}^d$ induces a map
$X(A, 0) \stackrel{\bar{\mu}}\to \mathbf{C}^d$ with $\bar{\mu}^{-1}(0) = Y(A, 0)$ and we have a commutative diagram
 
\begin{equation} 
\begin{CD} 
Y(A, \alpha) @>>> X(A, \alpha) \\ 
@VVV @VVV \\ 
Y(A, 0) @>>> X(A, 0) \\
@VVV @V{\bar{\mu}}VV  \\
0 @>>> \mathbf{C}^d     
\end{CD} 
\end{equation}  

Let $\{\mathbf{a}_1, ..., \mathbf{a}_n\}$ be the set of column vectors of $A$. Then its subset generates a subspace of $M_{\mathbf R}$. 
Consider all such subspaces $H_1$, ..., $H_l$ of $M_{\mathbf R}$ with {\em codimension $1$}. 
Then $\alpha \in M$ is called {\em generic} if $\alpha$ is not contained in any $H_i$.  
When $\alpha$ is generic, every point $p$ of $(\mu^{-1}(0))^{\alpha-ss}$ is $\alpha$-stable, i.e. $(\mu^{-1}(0))^{\alpha-ss} = 
(\mu^{-1}(0))^{\alpha-s}$ (cf. [Ko 2, Proposition 3.6], [HS, Proposition 6.2]). In particular, the stabilizer group 
$T^d$ of $p$ is a finite. This means that  $Y(A, \alpha)$ admits only quotient singularties. The symplectic 2-form $\omega$ on $\mathbf{C}^{2n}$ reduced to a symplectic orbifold 2-form $\omega_{Y(A, \alpha)}$ on 
the orbifold $Y(A, \alpha)$. In particular, $\omega_{Y(A, \alpha)}$ restricts to a usual symplectic 2-form on 
$Y(A, \alpha)_{reg}$. Similarly, $\omega$ is reduced to a symplectic 2-form $\omega_{Y(A, 0)}$ on 
$Y(A, 0)_{reg}$. We have $\omega_{Y(A, \alpha)} = \nu_{Y, \alpha}^*\omega_{Y(A, 0)}$. Then $Y(A, \alpha)$ has symplectic singularities by [Be, Proposition 2.4]. Hence $Y(A, 0)$ also has symplectic singularities. As $Y(A, 0)$ has a natural $\mathbf{C}^*$-action induced from 
the scaling $\mathbf{C}^*$-action on $\mathbf{C}^{2n}$, $Y(A, 0)$ is a conical symplectic variety with $wt(\omega_{Y(A, 0)}) = 2$. As $Y(A, \alpha)$ has 
symplectic singularities, $\mathrm{Sing}(Y(A, \alpha))$ has even codimension in $Y(A, \alpha)$ by [Ka].  By [Na 4, Corollary 03], 
$Y(A, \alpha)$ has terminal singularities if and only if $\mathrm{Codim}_{Y(A, \alpha)}\mathrm{Sing}(Y(A, \alpha)) \geq 4$. But it may happen that $\mathrm{Codim}_{Y(A, \alpha)}\mathrm{Sing}(Y(A, \alpha)) = 2$. 
Here is a criterion when $\mathrm{Codim}_{Y(A, \alpha)}\mathrm{Sing}(Y(A, \alpha)) = 2$ for a generic $\alpha$. 
\vspace{0.2cm}

{\bf Proposition 2}. {\em The following conditions are equivalent:}

(1) {\em For every generic element} $\alpha \in M$, $\mathrm{Codim}_{Y(A, \alpha)}\mathrm{Sing}(Y(A, \alpha)) = 2$.

(2) {\em For some $j_0 \in \{1, 2, ..., n\}$, the matrix $\bar{A} := (\mathbf{a}_1, ..., \mathbf{a}_{j_0 -1}, \mathbf{a}_{j_0 + 1}, ..., \mathbf{a}_n)$ 
satisfies the following conditions}:

(2-a): $\mathrm{rank}(\bar{A}) = d$, 

(2-b): {\em The homomorphism $\mathbf{Z}^{n-1} \stackrel{\bar{A}}\to \mathbf{Z}^d$ is not surjective.}
\vspace{0.2cm}

{\em Proof}. (1) $\Rightarrow$ (2): For a subset $J$ of $\{1,2, ..., n\}$, collect the column vectors $\mathbf{a}_j$ with $j \in J$ 
and form a matrix $d \times \vert J \vert$-matrix $A_J$. Take $J$ so that it is maximal among those with the following 
properties 

(i) $\mathrm{rank}(A_J) = d$

(ii) the homomorphism $\mathbf{Z}^{\vert J \vert} \stackrel{A_J}\to \mathbf{Z}^d$ is not surjective. 

We put $$(\mu^{-1}(0))_J := \{\mathbf{x} \in (\mu^{-1}(0))\: \vert \: z_j(\mathbf{x}) = w_j(\mathbf{x}) = 0, \:\: 
\forall j \notin J\}.$$ and define $$(\mu^{-1}(0))^{\alpha-s}_J := (\mu^{-1}(0))_J \cap (\mu^{-1}(0))^{\alpha-s}$$ 
  
Let us consider the quotient map $$\pi: (\mu^{-1}(0))^{\alpha-s} \to Y(A, \alpha).$$
We then have 
$$\mathrm{Sing}\: Y(A, \alpha) = \bigcup_J \pi ((\mu^{-1}(0))^{\alpha-s}_J).$$ 

{\bf Claim 3}. {\em Assume that $(\mu^{-1}(0))^{\alpha-s}_J \ne \emptyset$. Then we have:} 
\begin{itemize} 
\item $\dim \: (\mu^{-1}(0))^{\alpha-s}_J = 2\vert J \vert - d$
\item $\dim \: \pi ((\mu^{-1}(0))^{\alpha-s}_J ) =  2\vert J \vert - 2d$
\end{itemize} 

{\em Proof}. Let us consider the linear subspace $L$ of $\mathbf{C}^{2n}$ defined by $z_j = w_j = 0$, $j \notin J$. Then 
$L \cong \mathbf{C}^{2\vert J \vert}$ and we can take $z_j, w_j$ $(j \in J)$ as coordinates of $\mathbf{C}^{2\vert J \vert}$. 
Restrict the moment map $\mu$ to $L = \mathbf{C}^{2\vert J \vert}$. Then $\mu\vert_L: \mathbf{C}^{2\vert J \vert} \to \mathbf{C}^d$ 
is given by $$(\{z_j\}, \{w_j\})_{j \in J} \to \sum_{j \in J}\mathbf{a}_jz_jw_j.$$ Since $\mathrm{rank}\: A_J = d$, 
we can apply Lemma 1 to see that $(\mu\vert_L)^{-1}(0)$ 
is of complete intersection with $\dim = 2\vert J \vert - d$. (In the proof of Lemma 1, we only use the fact that $\mathrm{rank}\: A = d$.) 
Since $(\mu\vert_L)^{-1}(0) = (\mu^{-1}(0))_J$, we have $\dim \: (\mu^{-1}(0))^{\alpha-s}_J  = \dim \: (\mu^{-1}(0))_J = 
2\vert J \vert - d$. The second equality follows from the fact that every point $\mathbf{x} \in (\mu^{-1}(0))^{\alpha-s}_J$ has a finite stabilizer subgroup of $T^d$.  $\square$  \vspace{0.2cm}
  
Return to the proof (1) $\Rightarrow$ (2). Assume that $\mathrm{Codim}_{Y(A, \alpha)} \mathrm{Sing}\: Y(A, \alpha) = 2$ for a generic 
$\alpha$. Then, by Claim 3, there is a $J$ such that $\vert J \vert = n-1$ and $(\mu^{-1}(0))^{\alpha-s}_J \ne \emptyset$.  
If we put $j_0 = \{1,2, ..., n\} - J$, then this $j_0$ satisfies the condition of (2).  

(2) $\Rightarrow$ (1): We may assume $j_0 = n$ without loss of generality. Restrict the moment map $\mu$ to the $2n-2$-dimensional subspace $\mathbf{C}^{2n-2} := \{(z_1, ..., z_{n-1}, 0, w_1, ..., w_{n-1}, 0)\}$ of $\mathbf{C}^{2n}$. 
Then $$\mu\vert_{\mathbf{C}^{2n-2}}: \mathbf{C}^{2n-2} \to \mathbf{C}^d$$ is given by 
$$(z_1, ..., z_{n-1}, w_1, ..., w_{n-1}) \to \sum_{1 \le j  \le n-1}\mathbf{a}_jz_jw_j.$$  
Since $\mathrm{rank}\: (\mathbf{a}_1, ..., \mathbf{a}_{n-1}) = d$, we have $\dim (\mu\vert_{\mathbf{C}^{2n-2}})^{-1}(0) = 2n - 2 -d$ by 
Lemma 1. In other words, $\dim (\mu^{-1}(0) \cap \{z_n = w_n = 0\}) = 2n - 2 -d$.   
Let $\alpha \in M$ be a generic element. To prove that $\mathrm{Codim}_{Y(A, \alpha)}\mathrm{Sing}(Y(A, \alpha)) = 2$, it suffices 
to show that $(\mu^{-1}(0))^{\alpha-s} \cap \{z_n = w_n = 0\}) \ne \emptyset.$ 
In fact, if $(\mu^{-1}(0))^{\alpha-s} \cap \{z_n = w_n = 0\}) \ne \emptyset$, then $(\mu^{-1}(0))^{\alpha-s} \cap \{z_n = w_n = 0\})$ 
is a non-empty open subset of $\mu^{-1}(0) \cap \{z_n = w_n = 0\}$, which has dimension $2n - 2 -d$. 
Since $\mathrm{rank}\: (\mathbf{a}_1, ..., \mathbf{a}_{n-1}) = d$, we can take a sufficiently large positive integer $N$ so that $N\alpha$ can 
be written as  $$N\alpha = d_1\mathbf{a}_1 + \cdot\cdot\cdot + d_{n-1}\mathbf{a}_{n-1} \:\: (d_i \in \mathbf{Z})$$
We put $$J := \{j \: \vert \: 1 \le j \le n-1, \: d_j > 0\}, \:\: J' := \{j \: \vert \: 1 \le j \le n-1, \: d_j < 0\},$$ 
$$K := \{j \: \vert \: 1 \le j \le n-1, \: d_j = 0\}.$$
Take $\mathbf{x} = (z_1, ..., z_{n-1}, 0, w_1, ..., w_{n-1}, 0)$ so that 
\begin{itemize}
\item for $j \in J$, $z_j \ne 0$ and $w_j = 0$
\item for $j \in J'$, $z_j = 0$ and $w_j \ne 0$
\item for $j \in K$, $z_j \ne0$ and $w_j = 0$
\end{itemize}
Since $z_jw_j = 0$ for all $j$ with $1 \le j \le n-1$, we have $\mathbf{x} \in \mu^{-1}(0)$. 
Moreover, if we put 
$$f := \prod_{j \in J}z_j^{d_j} \prod_{j \in J'}w_j^{-d_j},$$ then $f \in \mathbf{C}[z_1, ..., z_n, w_1, ..., w_n]_{N\alpha}$ 
and $f(\mathbf{x}) \ne 0$. This means that $\mathbf{x} \in \mu^{-1}(0)^{\alpha-ss}$. 
Since $\alpha$ is generic, this also means that $\mathbf{x} \in \mu^{-1}(0)^{\alpha-s}$. 
Therefore, $\mathbf{x} \in \mu^{-1}(0)^{\alpha-s} \cap \{z_n = w_n = 0\}$; hence, 
$(\mu^{-1}(0))^{\alpha-s} \cap \{z_n = w_n = 0\}) \ne \emptyset.$
$\square$ \vspace{0.2cm}

\S {\bf 2}. 
\vspace{0.2cm}

{\bf (2.1)} We study in details a matrix $A$ satisfying the condition (2) of Proposition 2.
Let $A$ and $B$ be integer valued matrices of size $d \times n$ and $n \times (n-d)$ with an exact 
sequence $$0 \to \mathbf{Z}^{n-d} \stackrel{B}\to \mathbf{Z}^n \stackrel{A}\to \mathbf{Z}^d \to 0.$$ 
Let us consider two conditions respectively on $A$ and $B$: 

($\sharp_A$): Writing $A = (\mathbf{a}_1, \mathbf{a}_2, ..., \mathbf{a}_n)$ with the column vectors, the $d \times (n-1)$-matrix  
$\bar{A} = (\mathbf{a}_2, ..., \mathbf{a}_n)$ satisfies 

(a) $\mathrm{rank}\: \bar{A} = d$, and 

(b) the homomorphism $\mathbf{Z}^{n-1} \stackrel{\bar A}\to \mathbf{Z}^d$ is not surjective. 
\vspace{0.2cm}

($\sharp_B$): Writing  $$B = \left(\begin{array}{ccccccc} 
\mathbf{b}_1\\
 \mathbf{b}_2\\ 
... \\ 
\mathbf{b}_n  
\end{array}\right) \:\: \mathrm{with } \:\: \mathrm{row} \:\: \mathrm{vectors} \:\: \{\mathbf{b_i}\},$$ the 1-st row vector $\mathbf{b}_1$ is not a primitive vector 
with $\mathbf{b}_1 \ne 0$. 
\vspace{0.2cm}

{\bf Proposition 4}. {\em $A$ satisfies ($\sharp_A$) if and only if $B$ satisfies ($\sharp_B$).}
\vspace{0.2cm}

Before proving Proposition 4, we prepare a technical lemma, which will be also used in the proof of Lemma 6.  
\vspace{0.2cm}

{\bf Lemma 5}. {\em Assume that $A$ satisfies ($\sharp_A$).  If necessary,  replacing $A$ by $PA$ with an invertible integer valued  
matrix $P$ of size $d \times d$, we may assume that $A$ has the following form}  
$$A = \left(\begin{array}{cccccccc} 
 a_{11} &ma_{12} & ... & ma_{1n} \\
 a_{21} & a_{22} & ... &  a_{2n}\\ 
 ... & ... & ... & ...  \\ 
 a_{d1} & a_{d2} & ... & a_{dn}   
\end{array}\right)$$
{\em Here $m$ and $a_{ij}$ are both integers, $m > 1$ and $GCD(a_{11}, m) = 1$.    
Moreover, we may assume that the matrix} 
$$\left(\begin{array}{cccccccc} 
 1 &a_{12} & ... & a_{1n} \\
 0 & a_{22} & ... &  a_{2n}\\ 
 ... & ... & ... & ...  \\ 
 0 & a_{d2} & ... & a_{dn}   
\end{array}\right)$$ {\em determines a surjection $\mathbf{Z}^n \to \mathbf{Z}^d$.}  
\vspace{0.2cm}

{\em Proof}. 
Since $\mathrm{rank}\: \bar{A} = d$, we can take an invertible $d \times d$-matrix $P$ and an invertible $(n-1) \times (n-1)$-matrix $Q$ 
so that 
$$P\bar{A}Q = \left(\begin{array}{ccccccc} 
m_1 &  &  &  & 0 & ... & 0\\
 & m_2 &  &  & 0 & ... & 0\\ 
 &  & ... &  & 0 & ... & 0\\ 
 & &  & m_d & 0 & ...&0  
\end{array}\right),$$
where all $m_i$ are positive integers. Moreover, $m_1 > 1$ and each $m_i$ is a divisor of $m_{i-1}$. 
Then we have
$$PA\left(\begin{array}{cc} 
1 & \mathbf{0}\\
\mathbf{0} & Q  
\end{array}\right) = 
P\cdot(\mathbf{a}_1, \bar{A})\cdot\left(\begin{array}{cc} 
1 & \mathbf{0}\\
\mathbf{0} & Q  
\end{array}\right) = 
\left(\begin{array}{cccccccc} 
a'_{11} & m_1 &  &  &  & 0 & ... & 0\\
a'_{21} &  & m_2 &  &  & 0 & ... & 0\\ 
... & &  & ... &  & 0 & ... & 0\\ 
a'_{d1} & & &  & m_d & 0 & ...&0  
\end{array}\right), $$
where $$\left(\begin{array}{cccccc}
a'_{11} \\
a'_{21} \\
... \\
a'_{d1}\end{array}\right) =
P\mathbf{a}_1.$$ 
By assumption, the homomorphism $\mathbf{Z}^n \to \mathbf{Z}^d$ determined by 
$$PA\left(\begin{array}{cc} 
1 & \mathbf{0}\\
\mathbf{0} & Q  
\end{array}\right)$$ is a surjection. It follows from this fact that $m_2 = \cdot\cdot\cdot = m_d = 1$ and 
$GCD(a'_{11}, m_1) = 1$. We next write  
$$Q^{-1} = \left(\begin{array}{cc} 
Q'\\
Q''  
\end{array}\right)$$ with a $d \times (n-1)$-matrix $Q'$ and an $(n-d-1) \times n-1$-matrix $Q''$.  
Then we get 
$$PA =  \left(\begin{array}{cccccccc} 
a'_{11} & m_1 &  &  &  & 0 & ... & 0\\
a'_{21} &  & 1 &  &  & 0 & ... & 0\\ 
... & &  & ... &  & 0 & ... & 0\\ 
a'_{d1} & & &  & 1 & 0 & ...&0  
\end{array}\right)\cdot \left(\begin{array}{cc} 
1 & \mathbf{0}\\
\mathbf{0} & Q^{-1}  
\end{array}\right)$$
 
$$= \left(\begin{array}{cccccc}
a'_{11} \\
a'_{21} \\
... \\
a'_{d1}\end{array} 
\left(\begin{array}{cccccccc} 
m_1 &  &  &\\
 & 1 &  &  \\ 
 &  & ... &   \\ 
 & &  & 1   
\end{array}\right)Q' \right)$$ 
Here we can write
$$\left(\begin{array}{cccccccc} 
 m_1 &  &  &\\
  & 1 &  &  \\ 
  &  & ... &   \\ 
  & &  & 1   
\end{array}\right) Q'
= 
\left(\begin{array}{cccccccc} 
 m_1a'_{12} & ... & m_1a'_{1n} \\
 a'_{22} & ... &  a'_{2n}\\ 
 ... & ... & ...  \\ 
 a'_{d2} & ... & a'_{dn}   
\end{array}\right)$$ with suitable integers $a'_{ij}$. 
Then  
$$PA = \left(\begin{array}{cccccccc} 
 a'_{11} &m_1a'_{12} & ... & m_1a'_{1n} \\
 a'_{21} & a'_{22} & ... &  a'_{2n}\\ 
 ... & ... & ... & ...  \\ 
 a'_{d1} & a'_{d2} & ... & a'_{dn}   
\end{array}\right).$$
This calculation also shows that 
$$\left(\begin{array}{cccccccc} 
1 & m_1 &  &  &  & 0 & ... & 0\\
0 &  & 1 &  &  & 0 & ... & 0\\ 
... & &  & ... &  & 0 & ... & 0\\ 
0 & & &  & 1 & 0 & ...&0  
\end{array}\right)\cdot \left(\begin{array}{cc} 
1 & \mathbf{0}\\
\mathbf{0} & Q^{-1}  
\end{array}\right)
= 
\left(\begin{array}{cccccccc} 
 1 &m_1a'_{12} & ... & m_1a'_{1n} \\
 0 & a'_{22} & ... &  a'_{2n}\\ 
 ... & ... & ... & ...  \\ 
 0 & a'_{d2} & ... & a'_{dn}   
\end{array}\right)$$ 
Since 
$$\left(\begin{array}{cccccccc} 
1 & m_1 &  &  &  & 0 & ... & 0\\
0 &  & 1 &  &  & 0 & ... & 0\\ 
... & &  & ... &  & 0 & ... & 0\\ 
0 & & &  & 1 & 0 & ...&0  
\end{array}\right)$$ 
determines a surjection $\mathbf{Z}^n \to \mathbf{Z}^d$, 
$$\left(\begin{array}{cccccccc} 
 1 &m_1a'_{12} & ... & m_1a'_{1n} \\
 0 & a'_{22} & ... &  a'_{2n}\\ 
 ... & ... & ... & ...  \\ 
 0 & a'_{d2} & ... & a'_{dn}   
\end{array}\right)$$ also determines a surjection $\mathbf{Z}^n \to \mathbf{Z}^d$. 
By elementary transformations of column vectors, this matrix can be transformed to 
the matrix 
$$\left(\begin{array}{cccccccc} 
 1 & a'_{12} & ... & a'_{1n} \\
 0 & a'_{22} & ... &  a'_{2n}\\ 
 ... & ... & ... & ...  \\ 
 0 & a'_{d2} & ... & a'_{dn}   
\end{array}\right).$$
Hence this matrix also determines a surjection $\mathbf{Z}^n \to \mathbf{Z}^d$. 
$\square$   
\vspace{0.2cm}

{\em Proof of Proposition 4}.
Assume that $A$ satisfies ($\sharp_A$).  
We may assume that $A$ has the form in Lemma 5. 
If the 1-st row vector $\mathbf{b}_1$ of $B$ is zero, then $$AB = \bar{A} 
\left(\begin{array}{cccccccc} 
 \mathbf{b}_2 \\
 ...\\ 
 ...\\ 
 \mathbf{b}_n   
\end{array}\right) = 0$$. Since $\mathrm{rank}(\bar{A}) = d$,
$$\mathrm{we} \:\: \mathrm{have} \:\: \mathrm{rank}\left(\begin{array}{cccccccc} 
 \mathbf{b}_2 \\
 ...\\ 
 ...\\ 
 \mathbf{b}_n   
\end{array}\right) = n-1-d, \:\: \mathrm{which} \:\: \mathrm{means} \:\: \mathrm{that} \:\:  
\mathrm{rank}(B) = n-1-d.$$ This is a contradiction. Hence $\mathbf{b}_1 \ne 0$. 
The first row vector of $AB$ is then 
$$(a_{11}b_{11} + m(a_{12}b_{21} + \cdot\cdot\cdot ), a_{11}b_{12} + m(a_{12}b_{22} + \cdot\cdot\cdot ), ..., a_{11}b_{1,n-d} +m(a_{12}b_{2, n-d}  + \cdot\cdot\cdot ))$$ which 
must be zero. Since $GCD(a_{11}, m) = 1$, $m$ must divide $b_{1j}$ for any $j$. Therefore $B$ satisfies ($\sharp_B$). 

We prove the converse implication. Assume that $B$ satisfies ($\sharp_B$). 
Let us prove that $A$ satisfies (a) of ($\sharp_A$). Suppose to the contrary that 
$\mathrm{rank}(\bar{A}) \leq d-1$. Since $\mathrm{rank}(A) = d$, 
$\mathrm{rank}(\bar{A}) = d-1$. Then $$\mathrm{Ker}[\mathbf{Z}^{n-1} \stackrel{\bar{A}}\to \mathbf{Z}^d]$$ is a 
free abelian group of rank $n-d$. We take a generator 
$$\left(\begin{array}{cccccccc} 
 b'_{21} \\
 ...\\ 
 ...\\ 
 b'_{n1}   
\end{array}\right), \cdot\cdot\cdot,  
\left(\begin{array}{cccccccc} 
 b'_{2,n-d} \\
 ...\\ 
 ...\\ 
 b'_{n, n-d}   
\end{array}\right) \:\: \mathrm{of} \:\: \mathrm{this} \:\: \mathrm{group}, \mathrm{and} \:\: \mathrm{put}\:\:  
B' := \left(\begin{array}{cccccccc} 
 0 & ... & 0\\
 b'_{21} & ... & b'_{2, n-d}\\ 
 ... & ... & ...\\ 
 b'_{n1} & ... & b'_{n, n-d}   
\end{array}\right).$$ 
Then $\mathrm{rank}(B') = n-d$ and $AB' = 0$. This menas that $\mathrm{Im}(B)/\mathrm{Im}(B')$ is a 
torsion group. By the condition ($\sharp_B$), $b_{1j} \ne 0$ for some $j$. Then 
$$\left(\begin{array}{cccccccc} 
 b_{1, j} \\
 ...\\ 
 ...\\ 
 b_{n,j}   
\end{array}\right) \in \mathrm{Im}(B). \:\: \mathrm{Then} \:\: 
\mathrm{there} \:\:  \mathrm{is} \:\:  \mathrm{a} \:\: \mathrm{positive}\:\:  \mathrm{integer}  \:\: N \:\: 
\mathrm{such} \:\: \mathrm{that} \:\: 
N\left(\begin{array}{cccccccc} 
 b_{1, j} \\
 ...\\ 
 ...\\ 
 b_{n,j}   
\end{array}\right) \in \mathrm{Im}(B').$$ But this is absurd because $Nb_{1j} \ne 0$. 
We next prove that $A$ satisfies (b) of ($\sharp_A$). 
Assume, to the contrary, that $\bar{A}$ is a surjection. 
Consider the linear equation  
$$(\mathbf{a}_1, \bar{A})\left(\begin{array}{cccccccc} 
1 \\
x_2\\ 
 ...\\ 
x_n   
\end{array}\right) = \mathbf{a}_1 + \bar{A}
\left(\begin{array}{cccccccc} 
 x_2 \\
 ...\\ 
 ...\\ 
 x_n   
\end{array}\right) = 0$$
Since $\bar{A}$ is a surjection, this equation has a solution, say 
$$\left(\begin{array}{cccccccc} 
 c_2 \\
 ...\\ 
 ...\\ 
 c_n   
\end{array}\right). \:\: \mathrm{Then} \:\: 
\mathbf{c} = \left(\begin{array}{cccccccc} 
 1 \\
 c_2\\ 
 ...\\ 
 ... \\
 c_n   
\end{array}\right) \in \mathrm{Ker}(A).$$ By assumption, the first row vector $\mathbf{b}_1$ of $B$ is not primitive. 
This menas that, there is an integer $m$ with $m > 1$ such that, for any $\mathbf{c} \in \mathrm{Ker}(A)$, $m$ divides $c_1$. This is a contradiction. 
$\square$   \vspace{0.2cm}
  
{\bf (2.2)} For the matrix $A$ in Lemma 5, we introduce a new matrix $A'$ of size $(d+m-1) \times (n+m-1)$  
$$A' := \left(\begin{array}{ccccccccccc} 
 a_{11} & 0 & ... & ... & ... & 0 & a_{12} & a_{13} & ... & a_{1n} \\
 a_{21} & a_{21} & ... & ... & ... & a_{21} & a_{22} & a_{23} & ... & a_{2n}\\ 
 ...& ... & ...& ...& ...& ...& ...& ...& ... & ...\\ 
 a_{d1} & a_{d1} & ...& ... & ... & a_{d1} & a_{d2} & a_{d3} & ... & a_{dn}\\
1 & -1 & & & & & 0 & 0& ... & 0\\
 & 1 & -1 & & & & 0 & 0& ... & 0\\
 &  &  & ... &  &  & ... & ... & ... & ...\\
 & & & & 1& -1 & 0 & 0 & ... & 0 
\end{array}\right)$$
{\bf Lemma 6}. {\em $A'$ determines a surjection $\mathbf{Z}^{n+m-1} \to \mathbf{Z}^{d+m-1}$.}
   
{\em Proof}. By assumption the homomorphism $\mathbf{Z}^n \stackrel{A}\to \mathbf{Z}^d$ is a surjection. 
Hence, one can find a vector 
$$\left(\begin{array}{cccc} 
 c_1\\
 c_2\\ 
 ... \\
 c_n    
\end{array}\right) \:\:\:
\mathrm{such}\:\: \mathrm{that} \:\:\:   
\left(\begin{array}{cccccccc} 
 a_{11} &ma_{12} & ... & ma_{1n} \\
 a_{21} & a_{22} & ... &  a_{2n}\\ 
 ... & ... & ... & ...  \\ 
 a_{d1} & a_{d2} & ... & a_{dn}   
\end{array}\right)  
\left(\begin{array}{cccc} 
 c_1\\
 c_2\\ 
 ... \\
 c_n    
\end{array}\right)
= 
\left(\begin{array}{cccc} 
 1\\
 0\\ 
 ... \\
 0    
\end{array}\right).$$ 
Then we have 
$$\left(\begin{array}{ccccccccccc} 
 a_{11} & 0 & ... & ... & ... & 0 & a_{12} & a_{13} & ... & a_{1n} \\
 a_{21} & a_{21} & ... & ... & ... & a_{21} & a_{22} & a_{23} & ... & a_{2n}\\ 
 ...& ... & ...& ...& ...& ...& ...& ...& ... & ...\\ 
 a_{d1} & a_{d1} & ...& ... & ... & a_{d1} & a_{d2} & a_{d3} & ... & a_{dn}\\
1 & -1 & & & & & 0 & 0& ... & 0\\
 & 1 & -1 & & & & 0 & 0& ... & 0\\
 &  &  & ... &  &  & ... & ... & ... & ...\\
 & & & & 1& -1 & 0 & 0 & ... & 0 
\end{array}\right)
\left(\begin{array}{ccccccc} 
 c_1\\
 c_1\\ 
 ... \\
 c_1 \\
 mc_2 \\
 ... \\
 mc_n    
\end{array}\right)
= 
\left(\begin{array}{ccccccc} 
 1\\
 0\\ 
 ... \\
 ... \\
... \\
 ... \\
 0    
\end{array}\right).$$ 
We next take a vector
$$\left(\begin{array}{cccc} 
 d_1\\
 d_2\\ 
 ... \\
 d_n    
\end{array}\right) \:\:\:
\mathrm{such} \:\:  \mathrm{that} \:\:\:  
\left(\begin{array}{cccccccc} 
 a_{11} &ma_{12} & ... & ma_{1n} \\
 a_{21} & a_{22} & ... &  a_{2n}\\ 
 ... & ... & ... & ...  \\ 
 a_{d1} & a_{d2} & ... & a_{dn}   
\end{array}\right)  
\left(\begin{array}{cccc} 
 d_1\\
 d_2\\ 
 ... \\
 d_n    
\end{array}\right)
= 
\left(\begin{array}{cccc} 
 0\\
 1\\ 
 ... \\
 0    
\end{array}\right).$$ 
Then we have 
$$\left(\begin{array}{ccccccccccc} 
 a_{11} & 0 & ... & ... & ... & 0 & a_{12} & a_{13} & ... & a_{1n} \\
 a_{21} & a_{21} & ... & ... & ... & a_{21} & a_{22} & a_{23} & ... & a_{2n}\\ 
 ...& ... & ...& ...& ...& ...& ...& ...& ... & ...\\ 
 a_{d1} & a_{d1} & ...& ... & ... & a_{d1} & a_{d2} & a_{d3} & ... & a_{dn}\\
1 & -1 & & & & & 0 & 0& ... & 0\\
 & 1 & -1 & & & & 0 & 0& ... & 0\\
 &  &  & ... &  &  & ... & ... & ... & ...\\
 & & & & 1& -1 & 0 & 0 & ... & 0 
\end{array}\right)
\left(\begin{array}{ccccccc} 
 d_1\\
 0\\ 
 ... \\
 0 \\
 d_2 \\
 ... \\
 ... \\
 d_n    
\end{array}\right)
= 
\left(\begin{array}{cccccccccc} 
 * \\
 1 \\ 
 0 \\
 ... \\
 0 \\
 d_1 \\
 0 \\
 ... \\
 0    
\end{array}\right).$$ 

$$\mathrm{Since} \:\: \left(\begin{array}{ccccccc} 
 1\\
 0\\ 
 ... \\
 ... \\
 ... \\
 ... \\
... \\
 ... \\
 0    
\end{array}\right),  
\left(\begin{array}{cccccccccc} 
 * \\
 1 \\ 
 0 \\
 ... \\
 0 \\
 d_1 \\
 0 \\
 ... \\
 0    
\end{array}\right) \in \mathrm{Im}(A'), \:\: \mathrm{we} \:\: \mathrm{see} \:\: \mathrm{that} \:\: 
 \left(\begin{array}{cccccccccc} 
 0 \\
 1 \\ 
 0 \\
 ... \\
 0 \\
 d_1 \\
 0 \\
 ... \\
 0    
\end{array}\right) \in \mathrm{Im}(A').$$ 
$$\mathrm{Similarly}, \:\: \mathrm{we} \:\: \mathrm{see} \:\: \mathrm{that} \:\: 
\mathrm{Im}(A') \:\: \mathrm{contains} \:\: \mathrm{the} \:\: \mathrm{vectors} \:\: \mathrm{of} \:\: 
\mathrm{the} \:\: \mathrm{form} \:\: 
\left(\begin{array}{cccccccccc} 
 0 \\
 0 \\ 
 1 \\
 ... \\
 0 \\
 * \\
 0 \\
 ... \\
 0    
\end{array}\right), \cdot\cdot\cdot,  
\left(\begin{array}{cccccccccc} 
 0 \\
 0 \\ 
 ... \\
 0 \\
 1 \\
 * \\
 0 \\
 ... \\
 0    
\end{array}\right).$$ 
Recall that $$\left(\begin{array}{cccccccc} 
 1 &a_{12} & ... & a_{1n} \\
 0 & a_{22} & ... &  a_{2n}\\ 
 ... & ... & ... & ...  \\ 
 0 & a_{d2} & ... & a_{dn}   
\end{array}\right)$$ determines a surjection $\mathbf{Z}^n \to \mathbf{Z}^d$. Hence we can write 
$$\left(\begin{array}{cccccccccc} 
 0 \\
 a_{21} \\ 
 ... \\
 a_{d1}
\end{array}\right) = \lambda_1
\left(\begin{array}{cccccccccc} 
 1 \\
 0 \\ 
 ... \\
 0
\end{array}\right) + 
\lambda_2
\left(\begin{array}{cccccccccc} 
 a_{12} \\
 ... \\ 
 ... \\
 a_{d2}
\end{array}\right) + \cdot\cdot \cdot + \lambda_n
\left(\begin{array}{cccccccccc} 
 a_{1n} \\
 .... \\ 
 ... \\
 a_{dn}
\end{array}\right) \:\: \mathrm{with} \:\: \lambda_i \in \mathbf{Z}$$
Therefore we have    
$$\left(\begin{array}{cccccccccc} 
 0 \\
 a_{21} \\ 
 ... \\
 a_{d1}\\ 
 0 \\
 ... \\
 0
\end{array}\right) = \lambda_1
\left(\begin{array}{cccccccccc} 
 1 \\
 0 \\ 
 ... \\
 0 \\
 0 \\
 ... \\
 0 
\end{array}\right) + 
\lambda_2
\left(\begin{array}{cccccccccc} 
 a_{12} \\
 ... \\ 
 ... \\
 a_{d2} \\
 0 \\
 ... \\
 0
\end{array}\right) + \cdot\cdot \cdot + \lambda_n
\left(\begin{array}{cccccccccc} 
 a_{1n} \\
 .... \\ 
 ... \\
 a_{dn} \\
 0 \\ 
 ... \\ 
 0 
\end{array}\right).$$

$$\mathrm{Since} \:\: \left(\begin{array}{cccccccccc} 
 1 \\
 0 \\ 
 ... \\
 0 \\
 0 \\
 ... \\
 0 
\end{array}\right), 
\left(\begin{array}{cccccccccc} 
 a_{12} \\
 ... \\ 
 ... \\
 a_{d2} \\
 0 \\
 ... \\
 0
\end{array}\right), \cdot\cdot\cdot, 
\left(\begin{array}{cccccccccc} 
 a_{1n} \\
 .... \\ 
 ... \\
 a_{dn} \\
 0 \\ 
 ... \\ 
 0 
\end{array}\right) \in \mathrm{Im}(A'), \:\: \mathrm{we} \:\: \mathrm{see} \:\: \mathrm{that} \:\: 
\left(\begin{array}{cccccccccc} 
 0 \\
 a_{21} \\ 
 ... \\
 a_{d1}\\ 
 0 \\
 ... \\
 0
\end{array}\right) \in \mathrm{Im}(A').$$ 
$$\mathrm{Since} \:\:  
\left(\begin{array}{cccccccccc} 
 0 \\
 ... \\
 ... \\ 
 0 \\
 -1 \\ 
 1 \\
 0 \\
 ... \\
 0
\end{array}\right) = 
\left(\begin{array}{cccccccccc} 
 0 \\
 a_{21} \\ 
 ... \\
 a_{d1}\\ 
- 1 \\
1 \\
0 \\
 ... \\
 0
\end{array}\right) - 
\left(\begin{array}{cccccccccc} 
 0 \\
 a_{21} \\ 
 ... \\
 a_{d1}\\ 
0 \\
... \\
... \\
 ... \\
 0
\end{array}\right) \:\: \mathrm{and} \:\: 
\left(\begin{array}{cccccccccc} 
 0 \\
 a_{21} \\ 
 ... \\
 a_{d1}\\ 
- 1 \\
1 \\
0 \\
 ... \\
 0
\end{array}\right) \in \mathrm{Im}(A'), \:\:   
\left(\begin{array}{cccccccccc} 
 0 \\
 ... \\
 ... \\ 
 0 \\
 -1 \\ 
 1 \\
 0 \\
 ... \\
 0
\end{array}\right) \in \mathrm{Im}(A').$$
$$\mathrm{Similary} \:\: \mathrm{we} \:\: \mathrm{see} \:\: \mathrm{that} 
\left(\begin{array}{cccccccccc} 
 0 \\
 ... \\
 ... \\ 
 0 \\
 0 \\ 
 -1 \\
 1 \\
 0 \\
 0 \\
 ... \\
 0
\end{array}\right), 
\left(\begin{array}{cccccccccc} 
 0 \\
 ... \\
 ... \\ 
 0 \\
 0 \\ 
 0 \\
 -1 \\
 1 \\
 0 \\
 ... \\
 0
\end{array}\right), 
\cdot\cdot\cdot,  
\left(\begin{array}{cccccccccc} 
 0 \\
 ... \\
 ... \\ 
 0 \\
 0 \\ 
 0 \\
 0 \\
 0 \\
 0 \\
 ... \\
 -1
\end{array}\right) \in \mathrm{Im}(A').$$
Now let us consider the $(d+m-1) \times (d+m-1)$-matrix 
$$\left(\begin{array}{cccccccccccc}
1 & 0 & ... & ... & 0 & 0 & 0 & ... & ... & ... & 0 \\
0 & 1 & ... & ... & 0 & 0 & 0 & ... & ... & ... & 0 \\
... & 0 & ... & ... & ... & ... & ... & ... & ... & ... & ... \\
... & ... & ... & ... & ... & ... & ...  & ... & ... & ... & ... \\
... & 0 & ... & ... & 1 & 0  & 0 & ... & ... & ... & 0 \\
... & * & ... & ... & * & -1 & 0 & ... & ... & ... & 0 \\
... & 0 & ... & ... & 0 & 1 &  -1 & ... & ... &  ... & 0 \\
... & ... & ... & ... & ... & 0 &  1 & ... & ... &  ... & 0 \\
... & ... & ... & ... & ... & ... & 0 & ... & ... & ... & ... \\
... & ... & ... & ... & ... & ... & ... & ... & ... &  ... & 0 \\
0 & 0 & ... & ... & 0 &  0 & 0 & ... & ... &  ... & -1
\end{array}\right)$$ whose column vectors are all contained in $\mathrm{Im}(A')$. Then 
its determinant is $(-1)^{m-1}$, which implies the lemma. $\square$ 
\vspace{0.2cm}

By Proposition 4, we can take  
$$B = \left(\begin{array}{cccccccc} 
 mb_{11} &mb_{12} & ... & mb_{1,n-d} \\
 b_{21} & b_{22} & ... &  b_{2,n-d}\\ 
 ... & ... & ... & ...  \\ 
 b_{n1} & b_{n2} & ... & b_{n,n-d}   
\end{array}\right) \:\: \mathrm{with} \:\: m > 1.$$   
Now we put 
$$B' = \left(\begin{array}{cccccccc} 
 b_{11} & b_{12} & ... & b_{1,n-d} \\
 b_{11} & b_{12} & ... & b_{1,n-d} \\
 ... & ... & ... & ... \\
b_{11} & b_{12} & ... & b_{1,n-d} \\ 
 b_{21} & b_{22} & ... &  b_{2,n-d}\\ 
 ... & ... & ... & ...  \\ 
 b_{n1} & b_{n2} & ... & b_{n,n-d}   
\end{array}\right).$$
Here $B' $ is a $(n+m-1) \times (n - d)$-matrix and its first $m$ row vectors are 
$(b_{11}, b_{12}, ... , b_{1,n-d})$.  \vspace{0.2cm}

{\bf Lemma 7}. {\em The sequence $$0 \to \mathbf{Z}^{n-d} \stackrel{B'}\to \mathbf{Z}^{n+m-1} \stackrel{A'}\to 
\mathbf{Z}^{d+m-1} \to 0$$
is exact.}
\vspace{0.2cm}

{\em Proof}. One can directly check that $A'B' = 0$. It is enough to show that $(B')^t : \mathbf{Z}^{n+m-1} \to \mathbf{Z}^{n-d}$ 
is surjective. Since $B^t : \mathbf{Z}^n \to \mathbf{Z}^{n-d}$ is surjective, we see that $(B')^t$ is also surjective by the definition of 
$B'$. $\square$ \vspace{0.2cm} \newpage

\S {\bf 3}. 

Let $A$ and $B$ be the matrices in Proposition 4; namely, 
$$ A = \left(\begin{array}{cccccccc} 
 a_{11} &ma_{12} & ... & ma_{1n} \\
 a_{21} & a_{22} & ... &  a_{2n}\\ 
 ... & ... & ... & ...  \\ 
 a_{d1} & a_{d2} & ... & a_{dn}   
\end{array}\right), \:\:\: 
B = \left(\begin{array}{cccccccc} 
 mb_{11} &mb_{12} & ... & mb_{1,n-d} \\
 b_{21} & b_{22} & ... &  b_{2,n-d}\\ 
 ... & ... & ... & ...  \\ 
 b_{n1} & b_{n2} & ... & b_{n,n-d}   
\end{array}\right),$$ where $m$ is an integer with $m > 1$, and 
all $a_{ij}$, $b_{ij}$ are integers.     
Consider the affine space $\mathbf{C}^{2n}$ with coordinates $(z_m, ..., z_{n+m-1}, w_m, ..., w_{n+m-1})$ 
and put $$\omega := \sum_{m \le i \le n+m-1}dw_i \wedge dz_i.$$ As in \S 1, the matrix $A$ determines 
a $T^d$ action on $\mathbf{C}^{2n}$ (the coordinates suffixes being shifted) : 
$$(z_m, ..., z_{n+m-1}, w_m, ..., w_{n+m-1}) \to $$ $$(t_1^{a_{11}}t_2^{a_{21}}\cdot\cdot\cdot t_d^{a_{d1}}z_m, ..., t_1^{ma_{1n}}t_2^{a_{2n}}\cdot\cdot\cdot t_d^{a_{dn}}z_{n+m-1}, 
t_1^{-a_{11}}t_2^{-a_{21}}\cdot\cdot\cdot t_d^{-a_{d1}}w_m, ..., t_1^{-ma_{1n}}t_2^{-a_{2n}}\cdot\cdot\cdot t_d^{-a_{dn}}w_{n+m-1})$$
Let $$I = (i_{m+1}, ..., i_{n+m-1}) \in \mathbf{Z}^{n-1}_{\geq 0}, \:\: \mathrm{and}\:\: 
J := (j_{m+1}, ..., j_{n+m-1}) \in \mathbf{Z}^{n-1}_{\geq 0}$$ be $(n-1)$-tuples of non-negative integers such that 
$i_k = 0$ or $j_k = 0$ for every $k$ with $m+1 \le k \le n+m-1$. 
We then define $$(\mathbf{z}')^I(\mathbf{w}')^J := z_{m+1}^{i_{m+1}}\cdot\cdot\cdot z_{n+m-1}^{i_{n+m-1}}   
w_{m+1}^{j_{m+1}}\cdot\cdot\cdot w_{n+m-1}^{j_{n+m-1}}.$$

{\bf Claim 8}. {\em We have an inclusion} $$\mathbf{C}[z_m, ..., z_{n+m-1}, w_m, ..., w_{n+m-1}]^{T^d} \subset 
\mathbf{C}[z_m^m, w_m^m, z_mw_m, ..., z_{n+m-1}w_{n+m-1}, \{(\mathbf{z}')^I(\mathbf{w}')^J \}_{I, J}].$$ 

{\em Proof}. 
For $(\lambda_1, ..., \lambda_{n-d}) \in \mathbf{Z}^{n-d}$, we put 
$$ \left(\begin{array}{cccccccc} 
\tau_m \\
... \\
... \\
\tau_{n+m-1}   
\end{array}\right) := 
\lambda_1 
\left(\begin{array}{cccccccc} 
mb_{11} \\
b_{21} \\
... \\
b_{n1}   
\end{array}\right) + \cdot\cdot\cdot + 
\lambda_{n-d} 
\left(\begin{array}{cccccccc} 
mb_{1,n-d} \\
b_{2,n-d} \\
... \\
b_{n,n-d}   
\end{array}\right).$$ Then the Laurent monomial 
$z_m^{\tau_m}\cdot\cdot\cdot z_{n+m-1}^{\tau_{n+m-1}}$ is $T^d$-invariant. 
When $z_j^{\tau_j}$ with $\tau_j < 0$ appears in the Laurent monomial, substitute $w_j^{-\tau_j}$ for it. 
Then we get a $T^d$-invariant (usual) monomial of $z_m, ..., z_{n+m-1}, w_m, ..., w_{n+m-1}$. 
Such a monomial obtained in this way has the form 
$$(z_m^m)^i(\mathbf{z}')^I(\mathbf{w}')^J \:\: \mathrm{or} \:\: (w_m^m)^j(\mathbf{z}')^I(\mathbf{w}')^J. \:\:\: (*) $$
Conversely, the invariant ring $\mathbf{C}[z_m, ..., z_{n+m-1}, w_m, ..., w_{n+m-1}]^{T^d}$ is generated by such monomials and  
$z_mw_m$, ..., $z_{n+m-1}w_{n+m-1}$. Hence the claim holds. $\square$ 
\vspace{0.2cm}

Let $A'$ and $B'$ be the matrices introduced in (2.2): 
$$A' = \left(\begin{array}{ccccccccccc} 
 a_{11} & 0 & ... & ... & ... & 0 & a_{12} & a_{13} & ... & a_{1n} \\
 a_{21} & a_{21} & ... & ... & ... & a_{21} & a_{22} & a_{23} & ... & a_{2n}\\ 
 ...& ... & ...& ...& ...& ...& ...& ...& ... & ...\\ 
 a_{d1} & a_{d1} & ...& ... & ... & a_{d1} & a_{d2} & a_{d3} & ... & a_{dn}\\
1 & -1 & & & & & 0 & 0& ... & 0\\
 & 1 & -1 & & & & 0 & 0& ... & 0\\
 &  &  & ... &  &  & ... & ... & ... & ...\\
 & & & & 1& -1 & 0 & 0 & ... & 0 
\end{array}\right), \:\:\:\:   
B' = \left(\begin{array}{cccccccc} 
 b_{11} & b_{12} & ... & b_{1,n-d} \\
 b_{11} & b_{12} & ... & b_{1,n-d} \\
 ... & ... & ... & ... \\
b_{11} & b_{12} & ... & b_{1,n-d} \\ 
 b_{21} & b_{22} & ... &  b_{2,n-d}\\ 
 ... & ... & ... & ...  \\ 
 b_{n1} & b_{n2} & ... & b_{n,n-d}   
\end{array}\right).$$

Consider the affine space $\mathbf{C}^{2(n+m-1)}$ with coordinates $(z_1, ..., z_{n+m-1}, w_1 ..., w_{n+m-1})$ 
and put $$\omega' := \sum_{1 \le i \le n+m-1}dw_i \wedge dz_i.$$ As in \S 1, the matrix $A'$ determines 
a $T^{d + m-1}$ action on $\mathbf{C}^{2(n+m-1)}$. Then the moment map $$\mu': \mathbf{C}^{2(n+m-1)} \to \mathbf{C}^{d+m-1}$$ 
is given by 
$$\left(\begin{array}{cccccccccc}
z_1 \\
... \\ 
 ... \\
z_{n+m-1}\\ 
w_1 \\
... \\ 
... \\
 ... \\
w_{n+m-1}
\end{array}\right) \to
\left(\begin{array}{cccccccccc} 
 a_{11}z_1w_1 \:\:\:\:  \:\:\:\:  \:\:\:\:  \:\:\:\: \:\:\:\:  \:\:\:\:  \:\:\:\: \:\:\:\:  \:\:
+  a_{12}z_{m+1}w_{m+1} + \cdot\cdot\cdot + a_{1n}z_{n+m-1}w_{n+m-1} \\
 a_{21}z_1w_1  + \cdot\cdot\cdot + a_{21}z_mw_m + a_{22}z_{m+1}w_{m+1} + \cdot\cdot\cdot + a_{2n}z_{n+m-1}w_{n+m-1}  \\ 
 ... \\
 a_{d1}z_1w_1  +  \cdot\cdot\cdot + a_{d1}z_mw_m + a_{d2}z_{m+1}w_{m+1} + \cdot\cdot\cdot + a_{dn}z_{n+m-1}w_{n+m-1}\\ 
z_1w_1 - z_2w_2 \\
z_2w_2 - z_3w_3 \\
... \\
 ... \\
z_{m-1}w_{m-1} - z_mw_m
\end{array}\right)$$

Define a $d$-dimensional subspace $L$ of $\mathbf{C}^{d+m-1}$ by 
$$L := \{(u_1, ..., u_{d+m-1}) \in \mathbf{C}^{d+m-1} \:\vert \: u_{d+1} = \cdot\cdot\cdot = u_{d+m-1} = 0\}$$ and 
consider the inverse image $(\mu')^{-1}(L)$ of $L$ by $\mu'$. 
Since $(\mu')^*u_{d+1} = z_1w_1 - z_2w_2$, ..., $(\mu')^*u_{d+m-1} = z_{m-1}w_{m-1} - z_mw_m$, we have  
$$\mathbf{C}[(\mu')^{-1}(L)] = \mathbf{C}[z_1, ..., z_{n+m-1}, w_1, ... , w_{n+m-1}]/(z_1w_1 - z_2w_2, ..., z_{m-1}w_{m-1} - z_mw_m).$$  

{\bf Claim 9}. {\em The invariant ring
$\mathbf{C}[(\mu')^{-1}(L)]^{T^{d+m-1}}$ is a subring of} $$\mathbf{C}[z_1z_2\cdot\cdot\cdot z_m,\:\:\: w_1w_2\cdot\cdot\cdot w_m, 
\:\:\: z_1w_1, ..., 
z_{n+m-1}w_{n+m-1}, \{(\mathbf{z}')^I(\mathbf{w}')^J \}_{I, J}]/(\{z_iw_i - z_{i+1}w_{i+1}\}_{1 \le i \le m-1})$$

{\em Proof}. For $(\lambda_1, ..., \lambda_{n-d}) \in \mathbf{Z}^{n-d}$, we put 
$$ \left(\begin{array}{cccccccc} 
\tau_1 \\
... \\
... \\
\tau_{n+m-1}   
\end{array}\right) := 
\lambda_1 
\left(\begin{array}{cccccccc} 
b_{11} \\
... \\
b_{11} \\
b_{21} \\
... \\
b_{n1}   
\end{array}\right) + \cdot\cdot\cdot + 
\lambda_{n-d} 
\left(\begin{array}{cccccccc} 
b_{1,n-d} \\
... \\
b_{1,n-d} \\
b_{2,n-d} \\
... \\
b_{n,n-d}   
\end{array}\right).$$ Then the Laurent monomial 
$z_1^{\tau_1}\cdot\cdot\cdot z_{n+m-1}^{\tau_{n+m-1}}$ is $T^d$-invariant. 
When $z_j^{\tau_j}$ with $\tau_j < 0$ appears in the Laurent monomial, substitute $w_j^{-\tau_j}$ for it. 
Then we get a $T^d$-invariant (usual) monomial of $z_1, ..., z_{n+m-1}, w_1, ..., w_{n+m-1}$. 
Such a monomial obtained in this way has the form 
$$(z_1z_2\cdot\cdot z_m)^i(\mathbf{z}')^I(\mathbf{w}')^J \:\: \mathrm{or} \:\: (w_1\cdot\cdot w_m)^j(\mathbf{z}')^I(\mathbf{w}')^J.\:\:\:\: (**)$$
Conversely, the invariant ring $\mathbf{C}[(\mu')^{-1}(L)]^{T^{d+m-1}}$ is generated by such monomials and  
$z_1w_1$, ..., $z_{n+m-1}w_{n+m-1}$. Hence the claim holds. $\square$   
\vspace{0.2cm}

{\bf Remark 10}. The triplets $(i,  I, J)$, $(j, I, J)$  appeared in (**) coincide wth those appeared in (*) in the proof of Claim 8 because 
of the choice of $B$ and $B'$. 
\vspace{0.2cm}

Now let us consider a subring $R$ of $\mathbf{C}[u,v,w, z_{m+1}, ..., z_{n+m-1}, w_{m+1}, ..., w_{n+m-1}]/(uv - w^m)$ defined 
by $$R := \mathbf{C}[u, v, w, z_{m+1}w_{m+1}, ..., z_{n+m-1}w_{n+m-1}, \{(\mathbf{z}')^I(\mathbf{w}')^J \}_{I, J}]/(uv - w^m).$$

Then the ring  $$\mathbf{C}[z_1z_2\cdot\cdot\cdot z_m,\:\:\: w_1w_2\cdot\cdot\cdot w_m, 
\:\:\: z_1w_1, ..., 
z_{n+m-1}w_{n+m-1}, \{(\mathbf{z}')^I(\mathbf{w}')^J \}_{I, J}]/(\{z_iw_i - z_{i+1}w_{i+1}\}_{1 \le i \le m-1})
$$
appeared in Claim 9 is identified with $R$ by putting $u = m^mz_1z_2\cdot\cdot\cdot z_m$ and $v = w_1w_2\cdot\cdot\cdot w_m$ and 
$w = mz_mw_m$. 
 
On the other hand, the ring $$\mathbf{C}[z_m^m, w_m^m, z_mw_m, ..., z_{n+m-1}w_{n+m-1}, \{(\mathbf{z}')^I(\mathbf{w}')^J \}_{I, J}]$$ 
appeared in Claim 8 is also identified with $R$ by putting $u = z_m^m$, $v = w_m^m$ and $w = z_mw_m$.  

Therefore, these two rings are mutually identified via $R$. This identification induces an isomorphism   
$$\phi^*: \mathbf{C}[z_m, ..., z_{n+m-1}, w_m, ..., w_{n+m-1}]^{T^d} \to \mathbf{C}[(\mu')^{-1}(L)]^{T^{d+m-1}}$$ by Remark 10. 
Notice that $\phi^* (z_mw_m) = mz_mw_m$ and $\phi^*(z_jw_j) = z_jw_j$ for $j \geq m+1$. The isomorphism 
$\phi^*$ induces an isomorphism 
$$\phi: (\mu')^{-1}(L)/\hspace{-0.1cm}/_0T^{d+m-1} \to \mathbf{C}^{2n}/\hspace{-0.1cm}/_0T^d.$$

Let $\mu: \mathbf{C}^{2n} \to \mathbf{C}^d$ be the moment map for $(\mathbf{C}^{2n}, \omega)$. Then $\mu$ induces a map $\bar{\mu}: \mathbf{C}^{2n}/\hspace{-0.1cm}/_0T^d \to 
\mathbf{C}^d$, which is given by  
$$\overline{\left(\begin{array}{cccccccccc}
z_m \\
... \\ 
 ... \\
z_{n+m-1}\\ 
w_m \\
... \\ 
... \\
 ... \\
w_{n+m-1}
\end{array}\right)} \in  \mathbf{C}^{2n}/\hspace{-0.1cm}/_0T^d \to
\left(\begin{array}{cccccccccc} 
 a_{11}z_mw_m + ma_{12}z_{m+1}w_{m+1} + \cdot\cdot\cdot + ma_{1n}z_{n+m-1}w_{n+m-1} \\
 a_{21}z_mw_m + a_{22}z_{m+1}w_{m+1} + \cdot\cdot\cdot + a_{2n}z_{n+m-1}w_{n+m-1}  \\ 
 ... \\
 ... \\
 ... \\
 a_{d1}z_m w_m + a_{d2}z_{m+1}w_{m+1} + \cdot\cdot\cdot + a_{dn}z_{n+m-1}w_{n+m-1}  
\end{array}\right)$$ 
On the other hand, $\mu'$ induces a map $\bar{\mu'}_L: (\mu')^{-1}(L)/\hspace{-0.1cm}/_0T^{d+m-1} \to L$. Since $z_1w_1 = \cdot\cdot\cdot = w_mz_m$ 
on $(\mu')^{-1}(L)/\hspace{-0.1cm}/_0T^{d+m-1}$,  $\bar{\mu'}_L$ is given by 
$$\overline{\left(\begin{array}{cccccccccc}
z_1 \\
... \\ 
 ... \\
z_{n+m-1}\\ 
w_1 \\
... \\ 
... \\
 ... \\
w_{n+m-1}
\end{array}\right)} 
\to
\left(\begin{array}{cccccccccc} 
 a_{11}z_mw_m + a_{12}z_{m+1}w_{m+1} + \cdot\cdot\cdot + a_{1n}z_{n+m-1}w_{n+m-1} \\
 ma_{21}z_mw_m + a_{22}z_{m+1}w_{m+1} + \cdot\cdot\cdot + a_{2n}z_{n+m-1}w_{n+m-1}  \\ 
 ... \\
 ... \\
 ... \\
 ma_{d1}z_m w_m + a_{d2}z_{m+1}w_{m+1} + \cdot\cdot\cdot + a_{dn}z_{n+m-1}w_{n+m-1}  
\end{array}\right)$$ 
Then we have a $\mathbf{C}^*$-equivariant commutative diagram 
\begin{equation}
\begin{CD}
(\mu')^{-1}(L)/\hspace{-0.1cm}/_0T^{d+m-1} @>{\phi}>> \mathbf{C}^{2n}/\hspace{-0.1cm}/_0T^d\\ 
@V{\bar{\mu'}_L}VV  @V{\bar{\mu}}VV \\ 
L @>{\bar{\phi}}>> \mathbf{C}^d, 
\end{CD}
\end{equation}
where $\bar{\phi}$ is defined by $(u_1, u_2, ..., u_d) \to (mu_1, u_2, ..., u_d)$.  The horizontal maps $\phi$ and $\bar{\phi}$ are both isomorphisms.   
There is a relative symplectic form $\bar{\omega'}$ (with respect to $\bar{\mu'}$) on $(\mu')^{-1}(L)/\hspace{-0.1cm}/_0T^{d+m-1}$. 
Then $(\phi^{-1})^*\bar{\omega'}$ is a relative symplectic 2-form with respect to $\bar{\mu}$. 
Restrict $(\phi^{-1})^*\bar{\omega'}$ to $\bar{\mu}^{-1}(0) = Y(A, 0)$. Then it is a symplectic form on $Y(A, 0)$, which 
is denoted by $\omega'_{Y(A,0)}$. There is a natural symplectic form $\omega_{Y(A,0)}$ on $Y(A,0)$ defined by $\omega$. 
Since both $\omega_{Y(A,0)}$ and $\omega'_{Y(A,0)}$ have weight $2$, by [Na 1, Theorem 3.1], there is a $\mathbf{C}^*$-equivariant isomorphism 
of symplectic varieties
$$(Y(A, 0), \omega_{Y(A,0)}) \cong (Y(A, 0), \omega'_{Y(A,0)}).$$
In particular, we have $(Y(A, 0), \omega_{Y(A,0)}) \cong (\bar{\mu}^{-1}(0), \bar{\omega'}\vert_{\bar{\mu}^{-1}(0)})$ and 
$$\bar{\mu'}: \mathbf{C}^{2(m+n-1)}/\hspace{-0.1cm}/_0T^{d+m-1} \to \mathbf{C}^{d+m-1}$$ is regarded as a Poisson deformation of $(Y(A, 0), \omega_{Y(A,0)})$. 

Our strategy is as follows. Assume that $\mathrm{Codim}_{Y(A, \alpha)}\mathrm{Sing}(Y(A, \alpha) = 2$ for a generic $\alpha$. 
Then we take $A_1 := A'$ and consider $Y(A_1, 0)$. Then, as we have seen above, $Y(A_1, 0) \cong Y(A, 0)$. If $Y(A_1, \alpha_1)$ has only quotient terminal singularities for a generic $\alpha_1$, then the crepant partial resolution $Y(A_1, \alpha_1) \to Y(A_1, 0)$ gives a {\bf Q}-factorial terminalization of $Y(A, 0)$. If $\mathrm{Codim}_{Y(A_1, \alpha_1)}\mathrm{Sing}(Y(A_1, \alpha_1)) = 2$, then we put 
$A_2 := (A_1)'$ and consider $Y(A_2, \alpha_2)$ with a generic $\alpha_2$, and so on. We claim that this operation eventually terminates  
and we finally get $Y(A_k, \alpha_k)$ with only quotient terminal singularities. When we take $A_0$, $A_1$, $A_2$, ... starting with $A_0 := A$, 
we have exact sequences $$0 \to \mathbf{Z}^{n-d} \stackrel{B_i}\to \mathbf{Z}^{n+(m_1-1) + ... +(m_i -1)} \stackrel{A_i}\to 
\mathbf{Z}^{d + (m_1-1)+ ... + (m_i -1)} \to 0\:\:\: (i = 0,1, 2...).$$ What happens when we pass from $B_{i-1}$ to $B_i$ ?  As we have seen in Lemma 7, some row vector of $B_{i-1}$ 
has the form $(m_{i-1}b_1, ..., m_{i-1}b_{n-d})$. In $B_i$, it is broken into $m_{i-1}$ row vectors $(b_1, ..., b_{n-d})$. As long as 
there still exist a a non-primitive row in $B_i$, the operation does not terminate because of Proposition 4, and eventually   
all row vectors of $B_i$ become primitive at some stage, say $i = k$. Then $Y(A_k, \alpha_k)$ must have quotient 
terminal singularities again by Proposition 4. As a consequence, we have proved: \vspace{0.2cm}

{\bf Theorem 11}. {\em For a toric hyperk\"{a}hler variety $Y(A, 0)$, we can  take a suitable matrix $A^{\sharp}$ so that} \vspace{0.2cm}

(1) $(Y(A, 0), \omega_{Y(A,0)}) \cong (Y(A^{\sharp}, 0), \omega_{Y(A^{\sharp}, 0)})$ {\em as conical symplectic varieties, and}  \vspace{0.2cm}

(2) $Y(A^{\sharp}, \alpha^{\sharp}) \to Y(A^{\sharp}, 0)$ {\em is a {\bf Q}-factorial terminalization for a generic 
$\alpha^{\sharp}$.}
\vspace{0.2cm}

{\bf Example 12}. Let $$A = {\left(\begin{array}{cccccccccc}
1 & 0 & -2 & -2 \\
0 & 1 & -3 & -3
\end{array}\right)} \:\: \mathrm{and} \:\: 
B = {\left(\begin{array}{cccccccccc}
2 & 2 \\
3 & 3 \\
1 & 0 \\
0 & 1 
\end{array}\right)}.$$ These determine an exact sequence $$0 \to \mathbf{Z}^2 \stackrel{B}\to \mathbf{Z}^4 \stackrel{A}\to 
\mathbf{Z}^2 \to 0.$$
Since $A$ satisfies the conditions ($\sharp_A$) in Proposition 4, we have $\mathrm{Codim}_{Y(A, \alpha)}\mathrm{Sing}(Y(A, \alpha)) = 2$ for a generic $\alpha$. Accordng to (2.2), we introduce 
$$A' :=  {\left(\begin{array}{cccccccccc}
1 & 0 & 0 & -1 & -1 \\
0 & 0 & 1 & -3 & -3 \\
1 & -1 & 0 & 0 & 0
\end{array}\right)}.$$  
Permute the 1-st column and the 3-rd column of $A'$, and next permute the 1-st row and 
the 2-nd row of the resulting matrix. Then we get 
$$A_1 := {\left(\begin{array}{cccccccccc}
1 & 0 & 0 &  -3 & -3 \\
0 & 1 & 0 &  -1 & -1 \\ 
0  & 1 & -1 & 0 & 0
\end{array}\right)}.$$ Then $A_1$ together with 
$$B_1 := \left(\begin{array}{cccccccccc}
3 & 3 \\
1 & 1 \\ 
1 & 1 \\ 
1 & 0 \\
0 & 1
\end{array}\right)$$ determines an exact sequence 
$$0 \to \mathbf{Z}^2 \stackrel{B_1}\to \mathbf{Z}^5 \stackrel{A_1}\to \mathbf{Z}^3 \to 0.$$ 
The matrix  
${A}_1$ still satisfies the conditions (a) and (b), ($\sharp_A$) in Proposition 4. 
We put 
$$A_2 := \left(\begin{array}{cccccccccc}
1 & 0 & 0 & 0 & 0 & -1 & -1 \\
0 & 0 & 0 & 1 & 0 & -1 & -1 \\ 
0 & 0 & 0 & 1 & -1 & 0 & 0 \\ 
1 & -1 & 0 & 0 & 0 & 0 & 0 \\
0 & 1 & -1 & 0 & 0 & 0 & 0
\end{array}\right).$$ Then $A_2$ together with 
$$B_2 := \left(\begin{array}{cccccccccc}
1 & 1 \\
1 & 1 \\
1 & 1 \\
1 & 1 \\
1 & 1 \\ 
1 & 0 \\
0 & 1
\end{array}\right)$$ determines an exact sequence $$0 \to \mathbf{Z}^2 \stackrel{B_2}\to 
\mathbf{Z}^7 \stackrel{A_2}\to \mathbf{Z}^5 \to 0.$$
Since $A_2$ is unimodular, $Y(A_2, \alpha_2)$ is nonsingular for a generic $\alpha_2$. By construction,  
$Y(A, 0) \cong Y(A_1, 0) \cong Y(A_2, 0)$ and $Y(A_2, \alpha_2)$ gives a crepant resolution of $Y(A, 0)$. 
Notice that $\Sigma := \mathrm{Sing}(Y(A, 0))$ is irreducible and $\dim \Sigma = 2$. $Y(A, 0)$ has $A_4$-singularities along 
$\Sigma - \{\mathbf{0}\}$. $\square$  \vspace{0.2cm}

\S {\bf 4}.

Theorem 11 enables us to construct explicitly the universal Poisson deformation of $(Y(A, 0, ), \omega_{Y(A,0)})$ 
for  a toric hyperk\"{a}hler variety $Y(A, 0)$. For general properties of Poisson deformations, see [Na 2], [Na 3]. \vspace{0.2cm}

{\bf (4.1)}  First we treat the case when $Y(A, \alpha)$ gives a {\bf Q}-factorial terminalization of 
$Y(A, 0)$ for a generic $\alpha $. This is quite similar to the case where $A$ is unimodular and $Y(A, \alpha)$ is 
a crepant resolution of $Y(A, 0)$ (cf. [BLPW, 9.3], [Nag, Theorem 3.1]). The moment map $\mu: \mathbf{C}^{2n} \to \mathbf{C}^d$ induces a map 
$\bar{\mu}_{\alpha} : X(A, \alpha) \to \mathbf{C}^d$ with $\bar{\mu}_{\alpha}^{-1}(0) = Y(A, \alpha)$. $X(A, \alpha)$ is an orbifold, and 
$\bar{\mu}_{\alpha}$ is a family of orbifolds. There is a $\bar{\mu}_{\alpha}$-relative symplectic 2-form $\omega_{X(A, \alpha)}$ on 
the orbifold $X(A, \alpha)$, which restricts to a symplectic 2-form  $\omega_{Y(A, \alpha)}$ on the orbifold $Y(A, \alpha)$. 
In particular, $\omega_{X(A, \alpha)}$ determines a usual symplectic 2-form on the regualr 
part $X(A, \alpha)_{reg}$, which induces a Poisson structure $\{\:\:, \:\:\}_{X(A, \alpha)}$ on $X(A, \alpha)$. 
Similarly, $\omega_{X(A, \alpha)}$ determines a Poisson structure $\{\:\:, \:\:\}_{Y(A, \alpha)}$. The map $\bar{\mu}_{\alpha}$ 
is regarded as a Poisson deformation of $(Y(A, \alpha), \{\:\:, \:\:\}_{Y(A, \alpha)})$.     
\begin{equation} 
\begin{CD} 
Y(A, \alpha) @>>> X(A, \alpha) \\ 
@VVV @V{\bar{\mu}_{\alpha}}VV \\ 
0 @>>> \mathbf{C}^d = (\mathbf{t}^d)^*     
\end{CD} 
\end{equation}  
Each fber of $\bar{\mu}_{\alpha}$ is diffeomorphic (cf. [Ko 2, Proposition 3.6]), and $R^2(\bar{\mu}_{\alpha})_*\mathbf{C}$ is a constant sheaf of rank $b_2(Y(A, \alpha))$.  By using $\omega_{X(A, \alpha)}$, we get a period map for $\bar{\mu}_{\alpha}$: 
$$p: \mathbf{C}^d \to H^2(Y(A, \alpha), \mathbf{C}) \:\:\:\:  u \to [\omega_{X(A, \alpha)}\vert_{\bar{\mu}_{\alpha}^{-1}(u)}]$$
Let us introduce a map $\kappa$ from $(\mathbf{t}^d)^*$ to $H^2(Y(A, 0), \mathbf{C})$, called the Kirwan map.   
The quotient map $\mu^{-1}(0)^{\alpha-s} \to Y(A, 0)$ is an orbifold $T^d$-principal bundle. An 
element $\lambda \in \mathrm{Hom}_{alg.gp}(T^d, \mathbf{C}^*)$ determines an associated orbifold line bundle $L_{\lambda}$ 
on $Y(A, \alpha)$. We then have a homomorphism 
$$\mathrm{Hom}_{alg.gp}(T^d, \mathbf{C}^*) \to H^2(Y(A, \alpha),\mathbf{Q}), \:\:\: \lambda \to c_1(L_{\lambda}).$$ 
This map naturally extends to the map  from $(\mathbf{t}^d)^*$ to $H^2(Y(A, \alpha), \mathbf{C})$, which is the Kirwan map $\kappa$. 
The period map $p$ coincides with $\kappa$ (cf. [DH], [Ko 2, Proposition 6.1, (2)], [Lo, Proposition 3.2.1]). Moreover, the Kirwan map turns out to be an isomomorphism if any row of the matrix $B$ is 
nonzero (cf. [Ko 2, Proposition 6.1, (1)]). Therefore, $p$ is a linear isomorphism in our situation.

We next prove that $\bar{\mu}_{\alpha}$ is the universal Poisson deformation of $Y(A, \alpha)$.
Since $Y(A, 0)^{an}$ has only rational singularities and $Y(A, \alpha)^{an}$ is Stein, we have 
$$H^i(Y(A, \alpha)^{an}, \mathcal{O}_{Y(A, \alpha)^{an}}) = 0 \:\: \mathrm{for} \:\: i > 0.$$ Moreover, since $\mathrm{Codim}_{Y(A, \alpha)^{an}}\mathrm{Sing}(Y(A, \alpha)^{an}) \geq 4$ and $Y(A, \alpha)^{an}$ is Cohen-Macaulay, we have $H^i (Y(A, \alpha)^{an}_{reg}, \mathcal{O}_{Y(A, \alpha)^{an}_{reg}}) = 0$ for $i = 1, 2$ by the depth argument. 
By the exact sequences 
$$ 0 \to \mathbf{Z} \to \mathcal{O}_{Y(A, \alpha)^{an}} \to \mathcal{O}^*_{Y(A, \alpha)^{an}} \to 1, $$
$$ 0 \to \mathbf{Z} \to \mathcal{O}_{Y(A, \alpha)^{an}_{reg}} \to \mathcal{O}^*_{Y(A, \alpha)^{an}_{reg}} \to 1$$
we have isomorphisms $$\mathrm{Pic}(Y(A, \alpha)^{an}) \cong H^2(Y(A, \alpha), \mathbf{Z}) \:\: \mathrm{and} \:\:      
\mathrm{Pic}(Y(A, \alpha)^{an}_{reg}) \cong H^2(Y(A, \alpha)_{reg}, \mathbf{Z}).$$ On the other hand, we have 
$\mathrm{Pic}(Y(A, \alpha)^{an})\otimes_{\mathbf Z}\mathbf{Q} \cong 
\mathrm{Pic}(Y(A, \alpha)^{an}_{reg})\otimes_{\mathbf Z}\mathbf{Q}$ because $Y(A, \alpha)$ has only quotient singularities.  
Therefore, the restriction map $$H^2(Y(A, \alpha), \mathbf{Q}) \to H^2(Y(A, \alpha)_{reg}, \mathbf{Q})$$ is an isomorphism.        
Since $\mathrm{Codim}_{Y(A, \alpha)^{an}}\mathrm{Sing}(Y(A, \alpha)^{an}) \geq 4$, the 1-st order Poisson deformations of 
$Y(A, \alpha)$ are controlled by $H^2(Y(A, \alpha)_{reg}, \mathbf{C})$.  The Poisson deformation $\bar{\mu}_{\alpha}$ 
determines the Poisson Kodaira-Spencer map $$\tau: T_0(\mathbf{C}^d) \to H^2(Y(A, \alpha)_{reg}, \mathbf{C})(\cong 
H^2(Y(A, \alpha), \mathbf{C})),$$ which conicides with the period map $p$. In particular, $\tau$ is an isomorphism. 
This means that $\bar{\mu}_{\alpha}$ is the universal Poisson deformation of $Y(A, \alpha)$.  
\vspace{0.2cm}

{\bf (4.2)} As in (4.1) we assume that $Y(A, \alpha) \to Y(A, 0)$ is a {\bf Q}-factorial terminalization of $Y(A, 0)$ for a generic $\alpha$. 
Let us consider $B$. By Proposition 4, all row vectors of $B$ are primitive. We can assume that $B$ has the following form, after permuting the row vectors and changing their signs if necessary. 

$$B := \left(\begin{array}{cccccccccc}
\mathbf{b}_1 \\ 
... \\
\mathbf{b}_1 \\ 
\mathbf{b}_2 \\ 
... \\
\mathbf{b}_2\\
\mathbf{b}_3 \\
... \\
\mathbf{b}_3 \\ 
... \\
... \\
... \\ 
\mathbf{b}_r \\
... \\
\mathbf{b}_r
\end{array}\right),$$
where $\{\mathbf{b}_i\}$ are not mutually parallel, and each $\mathbf{b}_i$ appears in $d_i$ times.  
Then, as in [BLPW, 9.3] and [Nag, Theorem 3.11], the Weyl group $W$ of $Y(A, 0)$ is isomorphic to $\mathfrak{S}_{d_1} \times \cdot\cdot\cdot \times 
\mathfrak{S}_{d_r}$. The product $\mathfrak{S}_{d_1} \times \cdot\cdot\cdot \times 
\mathfrak{S}_{d_r}$ is naturally a subgroup of $\mathfrak{S}_{d_1+d_2 + ... +d_r}$ where $d_1 + ... + d_r = n$. 
Hence, an element of $\mathfrak{S}_{d_1} \times \cdot\cdot\cdot \times 
\mathfrak{S}_{d_r}$ is regarded as a permutation $\sigma$ of $\{1, 2, ..., n\}$. 
Then $\sigma \in W$ acts on $\mathbf{C}^{2n}$ by $$(z_1,..., z_n, w_1, ..., w_n) \to (z_{\sigma(1)}, ..., z_{\sigma(n)}, w_{\sigma(1)}, ..., w_{\sigma(n)})$$ 
The $W$ also acts on $X(A, 0) := \mathbf{C}^{2n}/\hspace{-0.1cm}/_0T^d$. The map $$\mathbf{C}^{2n} \to \mathbf{C}^n, \:\:\:\:  
(z_1, ..., z_n, w_1, ..., w_n) \to (z_1w_1, ..., z_nw_n)$$ descends to a map $f: X(A, 0) \to \mathbf{C}^n$. The map 
$\bar{\mu}$ factorizes as 
$$X(A, \alpha) \stackrel{f}\to \mathbf{C}^n \stackrel{A}\to \mathbf{C}^d$$
Introducing a $W$-action on $\mathbf{C}^n$ by natural permutations, the map $f$ is $W$-equivariant. As $W$ preserves $\mathrm{Ker}(A)$ by definition, the $W$-action on $\mathbf{C}^n$ descends to a 
$W$-action on $\mathbf{C}^d$.  In this way $\bar{\mu}$ is a $W$-equivariant map. 
We then have a commutative diagram   
\begin{equation} 
\begin{CD} 
X(A, 0) @>>> X(A, 0)/W \\ 
@V{\bar{\mu}}VV @V{\bar{\bar{\mu}}}VV \\ 
\mathbf{C}^d @>>> \mathbf{C}^d/W     
\end{CD} 
\end{equation}
The map  $\bar{\bar{\mu}}$ turns out to be the universal Poisson deformation of $Y(A, 0)$.   
\vspace{0.2cm} \newpage

{\bf(4.3)}

Let $Y(A, 0)$ be an arbitrary toric hyperk\"{a}hler variety. 
We can assume that $B$ has the following form, after 
permuting the row vectors and changing their signs if necessary. 
$$B := \left(\begin{array}{cccccccccc}
m_1\mathbf{b}_1 \\ 
... \\
m_{d_1}\mathbf{b}_1 \\ 
m_{d_1 + 1}\mathbf{b}_2 \\ 
... \\
m_{d_1 + d_2} \mathbf{b}_2\\
m_{d_1 + d_2 + 1}\mathbf{b}_3 \\
... \\
m_{d_1 + d_2 + d_3}\mathbf{b}_3 \\ 
... \\
... \\
... \\ 
m_{d_1 + ... + d_{r-1} + 1}\mathbf{b}_r \\
... \\
m_{d_1 + ... +d_{r-1} + d_r} \mathbf{b}_r
\end{array}\right),$$
where $m_i$ are positive integers, and $\mathbf{b}_i$ are primitive vectors, which are not mutually pararllel. 
Let $A^{\sharp}$ be the matrix in Theorem 11. Then $B^{\sharp}$ is the {\em primitivization} of $B$:  
$$B^{\sharp} := \left(\begin{array}{cccccccccc}
\mathbf{b}_1 \\ 
... \\
\mathbf{b}_1 \\ 
\mathbf{b}_2 \\ 
... \\
\mathbf{b}_2\\
\mathbf{b}_3 \\
... \\
\mathbf{b}_3 \\ 
... \\
... \\
... \\ 
\mathbf{b}_r \\
... \\
\mathbf{b}_r
\end{array}\right),$$
where each $\mathbf{b}_i$ appears in $m_{d_1 + ... + d_{i-1} + 1} + \cdot\cdot\cdot + 
m_{d_1 + ... +d_{i-1} + d_i}$ times.  
Therefore, the Weyl group $W$ of $Y(A^{\sharp}, 0)$ is isomorphic to 
$$\prod_{1 \le i \le r} \mathfrak{S}_{m_{d_1 + ... + d_{i-1} + 1} + \cdot\cdot\cdot + 
m_{d_1 + ... +d_{i-1} + d_i}}.$$  
Since $(Y(A, 0), \omega_{Y(A,0)}) \cong (Y(A^{\sharp}, 0), \omega_{Y(A^{\sharp}, 0)})$, 
$W$ is nothing but the Weyl group of $Y(A, 0)$. In (4.2) we have already constructed the universal Poisson 
deformation of $Y(A^{\sharp}, 0)$. Then it is also the universal Poisson deformation of $Y(A, 0)$. 
\vspace{0.2cm}

{\bf (4.4)} In the above, $Y(A^{\sharp}, \alpha^{\sharp})$ is nonsingular if and only if $A^{\sharp}$ is unimodular (cf. [HS, Proposition 6.2]). 
On the other hand, $A^{\sharp}$ is unimodular if and only if $B^{\sharp}$ is unimodular by the Gale duality.  
Therefore we have \vspace{0.2cm}

{\bf Corollary 13}. {\em A toric hyperk\"{a}hler variety $Y(A, 0)$ has a projective crepant resolution if and only 
if the primitivization $B^{\sharp}$ of $B$ is unimodular.}
\vspace{0.2cm}

{\em Proof}. Since we have already proved the ``if'' part, we only have to prove the ``only if'' part. 
Assume that $B^{\sharp}$ is not unimodular. Then $Y(A^{\sharp}, \alpha^{\sharp})$ is a singular projectve {\bf Q}-factorial 
terminalization for a generic $\alpha^{\sharp}$. This means that any projective {\bf Q}-factorial terminalization of $Y(A^{\sharp}, 0)$ is 
singular\footnote{More strongly, one can prove that any {\bf Q}-factorial terminalization of $Y(A^{\sharp}, 0)$ is obtained as the hyperk\"{a}hler  reduction $Y(A^{\sharp}, \alpha^{\sharp})$ with a generic  
$\alpha^{\sharp}$.} by [Na 2, Corollary 25].   
  $\square$ 
\vspace{0.5cm}

\begin{center}
{\bf References} 
\end{center}

[Be] Beauville, A.: Symplectic singularities, Invent. Math. {\bf 139} (2000), no.3, 541-549
\vspace{0.15cm}


 
[BD] Bielawski, R.,  Dancer, A.: 
The geometry and topology of toric hyperk\"{a}hler manifolds.
Comm. Anal. Geom. {\bf 8} (2000), no. 4, 727 - 760
\vspace{0.15cm}

[BLPW] Braden, T., Licata, A., Proudfoot, N., Webster, B.: 
Quantizations of conical symplectic resolutions II: category O and symplectic duality. 
with an appendix by I. Losev.
Ast\'{e}risque No. {\bf 384} (2016), 75 - 179.
\vspace{0.15cm}
 
[DS] Duistermaat, J. J., Heckman, G., J.:  On the variation in the cohomology of the symplectic form of the reduced phase space. 
Invent. Math. {\bf 69} (1982), no. 2, 259 - 268 
\vspace{0.15cm}


[Go] Goto, R.: On toric hyper-K\"{a}hler manifolds given by the hyper-K\"{a}hler quotient method, in Infinite Analysis, World Scientific, 1992, pp. 317 - 338
\vspace{0.15cm}


[HS] Hausel, T.; Sturmfels, B.: 
Toric hyperk\"{a}hler varieties. 
Doc. Math. {\bf 7} (2002), 495 - 534.
\vspace{0.15cm}
     

[Ka] Kaledin, D.: Symplectic singularities from the Poisson point of view.
J. Reine Angew. Math. {\bf 600} (2006), 135 - 156.
\vspace{0.15cm}

[Ko 1] Konno, H.: 
Cohomology rings of toric hyperk\"{a}hler manifolds.
Internat. J. Math. {\bf 11} (2000), no. 8, 1001 - 1026.
\vspace{0.15cm}

[Ko 2] Konno, H.: The geometry of toric hyperk\"{a}hler varieties. Toric topology, 241- 260,
Contemp. Math., {\bf 460}, Amer. Math. Soc., Providence, RI, 2008.
\vspace{0.15cm}

[Lo] Losev, I.: Isomorphisms of quantizations via quantization of resolutions. 
Adv. Math. {\bf 231} (2012), no. 3-4, 1216 - 1270. 
\vspace{0.15cm}

[Nag] Nagaoka, T.: The universal Poisson deformation of hypertoric varieties and some classification results. Pacific J. Math. {\bf 313} (2021), no. 2, 459 - 508.
\vspace{0.15cm}

[Na 1] Namikawa, Y.: Equivalence of symplectic singularities. 
Kyoto J. Math. {\bf 53} (2013), no. 2, 483 - 514. 
\vspace{0.15cm}

[Na 2] Namikawa, Y.: Flops and Poisson deformations of symplectic varieties. Publ. Res. Inst. Math. Sci. {\bf 44} (2008), no. 2, 259 - 314.
\vspace{0.15cm}

[Na 3] Namikawa, Y.: Poisson deformations of affine symplectic varieties. Duke Math. J. {\bf 156} (2011), no. 1, 51 - 85.
\vspace{0.15cm}
 
[Na 4] Namikawa, Y.: Birational geometry for the covering of a nilpotent orbit closure. 
Selecta Math. (N.S.) {\bf 28} (2022), no. 4, Paper No. 75, 59 pp.
\vspace{0.5cm}



 



\begin{center}
Research Institute for Mathematical Sciences, Kyoto University, Oiwake-cho, Kyoto, Japan

E-mail address: namikawa@kurims.kyoto-u.ac.jp  
\end{center}

\end{document}